# A hybrid probabilistic domain decomposition algorithm suited for very large-scale elliptic PDEs

Francisco Bernal*, Jorge Morón-Vidal†, and Juan A. Acebrón+

**Abstract**

State of the art domain decomposition algorithms for large-scale boundary value problems (with $M \gg 1$ degrees of freedom) suffer from bounded strong scalability because they involve the synchronisation and communication of workers inherent to iterative linear algebra. Here, we introduce PDDSparse, a different approach to scientific supercomputing which relies on a "Feynman-Kac formula for domain decomposition". Concretely, the interfacial values (only) are determined by a stochastic, highly sparse linear system $G(\omega)\vec{u} = \vec{b}(\omega)$ of size $O(\sqrt{M})$, whose coefficients are constructed with Monte Carlo simulations—hence embarrassingly in parallel. In addition to a wider scope for strong scalability in the deep supercomputing regime, PDDSparse has built-in fault tolerance and is ideally suited for GPUs. A proof of concept example with up to 1536 cores is discussed in detail.

* Department of Mathematics, Carlos III University of Madrid, Spain (franciscomanuel.bernal@uc3m.es)
† Department of Mathematics, Carlos III University of Madrid, Spain (jmoron@math.uc3m.es)
+ ISCTE-IUL and INESC-ID, Instituto Superior Técnico, Universidade de Lisboa, Portugal & Department of Mathematics, Carlos III University of Madrid, Spain (juan.acebron@iscte-iul.pt)

## 1 Introduction

**Background.** The need for numerical solutions of large-scale scale models based on bounday value problems (BVPs) is ubiquitous in applied mathematics. Owing to their sheer discretisation size (billions of degrees of freedom, or DoFs), not only do such BVPs call for supercomputers (i.e. massive distributed computers), but also for sophisticated algorithms designed to make the best use of the available computational resources. State of the art algorithms [13, 29, 30] consist in partitioning the original BVP domain into a tessellation of artificial subdomains which are small enough that they can be tackled by individual cores/threads of the supercomputer. However, because the subdomain-restricted BVPs are bound to be ill-posed—for lack of the correct subdomain boundary conditions (BCs)—an iterative loop is unavoidable, in which subdomain BCs are updated from their neighbours, until convergence. Such intercommunication is intrinsically sequential and—by virtue of Amdahl's law—sets an upper bound to the strong scalability potential of state of the art algorithms. In other words, the solution of a BVP with a given accuracy target can be accelerated only so much regardless of the available hardware.

A radical alternative was introduced with probabilistic domain decomposition (PDD) [1]. The insight is to solve for the artificial interfaces first, by resorting to the probabilistic representation of the BVP—epitomised by the Feynman-Kac formula [18]. Since the latter can be solved embarrassingly in parallel by means of Monte Carlo simulations, the subdomain-restricted BVPs are rendered well posed, and can be solved independently from one another.

Notwithstanding its potentially full strong scalability (which sets it apart from the state of the art), PDD poses some nontrivial issues of its own:

- The interfacial solution carries the statistical error of Monte Carlo and the weak error of the stochastic differential equation (SDE) solver. Both types of errors converge notoriously slowly; typically as the square root of the number of trajectories (in the first case), and linearly with the timestep (in the second) [22]. This leads to lengthy simulations.

- In order to integrate a Feynman-Kac trajectory until hitting the real boundary, the core in charge of it must know the BVP coefficients throughout the complete domain. When such coefficients are available as interpolatory (look-up) tables (instead of formulas), they must fit in each individual core's memory. Otherwise, whenever a Feynman-Kac trajectory strays out of the region "known" to the integrating core, the latter must hand over control of that trajectory to another core which holds the required information. This exchange introduces intercore communication—causing Amdahl's law to rear its ugly head.

**This paper's contribution.** Here, we present a novel algorithm partaking in the PDD paradigm, called *PDDSparse*, which both overcomes the "limited core memory" issue and dramatically speeds up the simulations. PDDSparse can be applied to linear elliptic or parabolic BVPs. However, in this paper, we shall address two-dimensional elliptic BVPs only, for the sake of simplicity.

PDDSparse stems from the observation that interfacial nodes (or *knots*) close to one another give rise to rather similar ensembles of trajectories, which (in PDD) must be integrated up to the real boundary. On the other hand, trajectories somewhat "lose memory" of their starting point once they drift far away from it. It would be desirable to split the contribution of stochastic trajectories into a "local part" (which concentrates the effect of the starting point), and a "distant part" (which can be reused by clusters of nearby knots).

For the purpose of illustration, consider the square domain $\Omega$ depicted in Figure 1. It has been partitioned into a grid of nonoverlapping square subdomains. The same partition is used by PDD (left), and PDDSparse (right).

- With PDD, Feynman-Kac trajectories from every *interfacial knot* are integrated up to their first-exit point on $\partial\Omega$, yielding the *nodal values* $u_1, \ldots, u_n$, where $u_k \approx u(\mathbf{x}_k)$ is the (approximate) solution on the $k^{th}$ knot.

- With PDDSparse, trajectories from $\mathbf{x}_k$ are integrated until hitting the boundary $\partial\Pi(\mathbf{x}_k)$ of the *patch* around $\mathbf{x}_k$, $\Pi(\mathbf{x}_k)$. The knot's patch is the union of the subdomains that the knot belongs to. The (unknown) BCs on $\partial\Pi(\mathbf{x}_k)$ can be approximated through linear interpolation of the knots on it. Indeed, it will be shown in Section 3 that $u_k \in span(u_1^{(k)}, u_2^{(k)}, \ldots, u_{n_k}^{(k)})$, where $u_j^{(k)}$ are the nodal values of the nodes on $\partial\Pi(\mathbf{x}_k)$. This is possible because the PDE is linear and elliptic and $\partial\Pi(\mathbf{x}_k)$ surrounds $\mathbf{x}_k$. Repeating the procedure for each one of the $n$ interfacial knots (for knots close to $\partial\Omega$ there is a slight variation), a linear system $G\vec{u} = \vec{b}$ arises, in which row $k$ of the *PDDSparse matrix G* stems from applying the Feynman-Kac formula along with linear interpolation on the patch of $\mathbf{x}_k$. It can, in fact, be regarded as a discrete Feynman-Kac formula for domain decomposition.

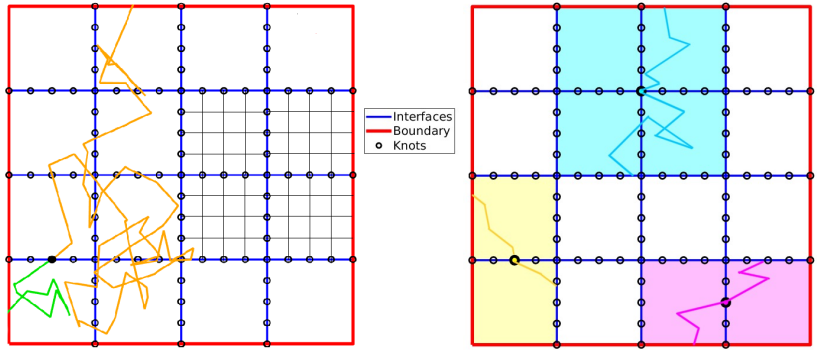

Figure 1: PDD (left) ignores the interior discretisation of subdomains and directly solves the BVP on the interfacial knots only, by means of stochastic trajectories like the orange and green ones. They must be integrated until hitting the boundary. With PDDSparse (right), trajectories are integrated only up to the patch boundary. (Three patches are highlighted in colour. Two trajectories each are shown from their respective central knots.)

Integrating within patches offers critical advantages:

- Since the patches are only up to four subdomains big, the "memory issue" of PDD is fixed. Handling trajectories is guaranteed to be embarrassingly parallel, regardless of the number of subdomains, $M_{subs}$.

- Trajectories take much less (by a factor or the order of $1/M_{subs}$) to hit their patch boundary than the actual one. Moreover, the weak error bias and the variance of the associated Feynman-Kac functionals typically decrease with the integration area. The three effects build up in such a way that the computational cost is dramatically lowered with respect to PDD.

- The discretisation setup of PDDSparse allows for tailored schemes which improve and/or control the errors. Variance reduction based on pathwise control variates [6] is particularly amenable to PDDSparse, as it is the estimation of the pointwise leading bias constants.

PDDSparse results in a linear system *for the nodal, or interfacial, solutions*. It must be emphasised that this system of size $n$ is orders of magnitude smaller than the one involving all of the original problem DoFs, $M$. As a conservative example, assume that it takes $1000 \times 1000$ DoFs to discretise one of the square subdomains in Figure 1 for the desired accuracy.[1] With standard domain decomposition, that subdomain would contribute one million DoFs to the global linear system (which involves all the knots that $\Omega$ is discretised into). On the other hand, it just contributes ca. 2000 DoFs (half of the boundary knots[2]) to the PDDSparse linear system. With tortuous domains, the difference is bigger, since there are less interior (or "floating") subdomains and hence proportionally fewer interfacial knots.

There are very successful state of the art algorithms which also lead to systems for the interfacial knots only (called the *skeleton* in that context). They are known as iterative substructuring, Schur complement, or nonoverlapping methods. (For instance, Balancing Neumann-Neumann [23], Balancing Domain Decomposition by Constraints [12], and FETI [17].) Unfortunately, the stiffness matrix loses its sparsity upon shrinkage, and involves the inverse of the interior operators. The explicit calculation of the Schur complement is expensive and calls for a large amount of memory since it is much denser than the original, $O(M)$, one. In order to avoid forming it, its action (the products of inverse matrices and vectors) are replaced by local linear systems [27] and solved via preconditioned Krylov iterations [28]. The choice of the preconditioner is critical and often problem-dependent, and requires solving delicate auxiliary problems on the interfaces [5, 10, 15] or on the overlaps at each iteration. Moreover, the number of iterations of the preconditioned system typically grows with the number of subdomains [30], which goes against the ideal strong scalability scenario of matching subdomains and cores. For that reason, coarse spaces are needed that accelerate convergence [14, 16, 23]. In practice, coarse space tasks become the sequential bottleneck *even for weak scalability*—limiting the number of cores which can be exploited strongly in parallel to a few tens of thousands, with very careful implementations [2] (see also [13, chapter 8]). In contrast, the PDDSparse system is an explicit, sparse, structured, rather general, "one shot", $O(\sqrt{M})$ system. Even though domain decomposition is an extremely fruitful and vibrant area, it is unlikely that incremental improvements to strong scalability can keep pace with the hardware revolution, at least in the short term. This justifies looking for alternative (or complementary) ideas which rely on entirely different formulations with untapped potential, such as PDDSparse.

---

[1] In that sketch, the fine grid consists of just $5 \times 5$ DoFs per subdomain.
[2] Most knots belong to two squares but they show up just once in the PDDSparse linear system.

Notwithstanding the attainable shrinkage, PDDSparse still results in a linear system—and that fact will eventually put an end to its strong scalability, as happens with the state of the art. However, this limit lies potentially deep into the supercomputing range (i.e. it is expected to occur with many more cores than now). Adding to this its natural fault tolerance and its suitability for GPUs, PDDSparse plays to the strengths of modern supercomputer architecture, namely GPU-accelerated cores with millions of cores altogether.[3]

The stability—and even the invertibility—of the PDDSparse matrix $G$ is a legitimate concern. Our numerical observations so far are distinctly encouraging. Theoretical investigations, to be reported elsewhere [9], indeed suggest that $G$ enjoys an underlying structure (beyond sparsity) which keeps the condition number growth with $n$ and $M_{subs}$ under control, and preconditionable.

**Organisation of the remainder of the paper.** The next section is a briefing on several ingredients of PDDSparse. Concretely, one-dimensional radial basis function interpolation; the Feynman-Kac formula required for Section 3; and an overview of the associated numerical methods for bounded SDEs. In Section 3, the new algorithm PDDSparse is formulated for 2D (linear) elliptic BVPs with Dirichlet BCs. (Previously, the discretisation and terminology are explained.) A comprehensive discussion of the features of PDDSparse follows. Section 4 is devoted to the more salient implementation aspects of PDDSparse, and lists the pseudocode. One composed numerical experiment is reported in Section 5 and criticised, paving the way for the conclusions.

## 2 Preliminaries

### 2.1 One-dimensional RBF interpolation

Let $J = [a, b] \subset \mathbb{R}$ be a closed interval and $u : J \mapsto \mathbb{R}$ a smooth function. Given the $m$ interpolatory data $\{(z_i, u_i)\}_{i=1}^m$, where $a = z_1 < z_2 < \ldots < z_m = b$ and $u_i = u(z_i)$; and a radial basis function (RBF) $\phi_c(r)$, the RBF interpolant of $u$ is

$$\tilde{u}(z) = \sum_{i=1}^m \alpha_i \phi_c(|z - z_i|), \tag{1}$$

---
[3]https://www.top500.org/lists/top500/2022/11/

where the real coefficients $\alpha_1, \ldots, \alpha_m$ are determined by the (interpolation) condition $\tilde{u}(z_i) = u_i$, i.e. by collocating (1) on every knot of the *pointset* $\{z_i\}$:

$$\begin{bmatrix} \alpha_1 \\ \alpha_2 \\ \vdots \\ \alpha_m \end{bmatrix} = \begin{bmatrix} \phi_c(0) & \phi_c(|z_1 - z_2|) & \cdots & \phi_c(|z_1 - z_m|) \\ \phi_c(|z_2 - z_1|) & \phi_c(0) & \cdots & \phi_c(|z_2 - z_m|) \\ \vdots & \vdots & \ddots & \vdots \\ \phi_c(|z_m - z_1|) & \phi_c(|z_m - z_2|) & \cdots & \phi_c(0) \end{bmatrix}^{-1} \begin{bmatrix} u_1 \\ u_2 \\ \vdots \\ u_m \end{bmatrix} =: \Phi^{-1} \begin{bmatrix} u_1 \\ u_2 \\ \vdots \\ u_m \end{bmatrix}.$$
(2)

The *RBF interpolation matrix* $\Phi$ above is symmetric, and invertible if $\phi_c$ is positive definite [32]. In this paper, we have used the inverse multiquadric,

$$\phi_c(r) = \frac{1}{\sqrt{r^2 - c^2}},$$
(3)

where the real parameter $c$ is called the *shape parameter*. Inserting $\alpha_1, \ldots, \alpha_m$ from (2) into (1) yields the *Lagrange form* of the RBF interpolant,

$$\tilde{u}(z) = \sum_{i=1}^{m} \left( \sum_{j=1}^{m} \Phi_{ij}^{-1} u_j \right) \phi_c(|z - z_i|),$$
(4)

where $\Phi_{ij}^{-1}$ stands for the $(i, j)$ entry of the square symmetric matrix $\Phi^{-1}$ in (2).

The interpolant in (4) can be rewritten as

$$\tilde{u}(z) = \sum_{j=1}^{m} u_j H_j(z), \qquad \text{where} \qquad H_j(z) = \sum_{i=1}^{m} \Phi_{ij}^{-1} \phi_c(|z - z_i|)$$
(5)

is the *cardinal function* associated to $z_j$. It bears that name because (as it is not difficult to show) $H_j(z_i) = \delta_{ij}$ (Kronecker's delta).

If $u(z) \in C^\infty$, one-dimensional RBF interpolants are spectrally convergent with respect to $1/m$ and to the shape parameter $c$ for many RBFs, such as (3). They do not require that the stencil $\{z_i\}$ be equispaced, thus retaining scope for adaptivity. (With an equispaced pointset, though, cardinal functions tend to the sinc function, which lays the groundwork for the theoretical study in [9].)

## 2.2 The Feynman-Kac formula

Let $\Omega \subset \mathbb{R}^2$ be an open bounded domain with boundary $\partial\Omega$. Consider the second order linear elliptic BVP with Dirichlet boundary conditions:

$$\begin{cases} \underbrace{\frac{1}{2}a_{xx}(\mathbf{x})\frac{\partial^2 u}{\partial x^2} + a_{xy}(\mathbf{x})\frac{\partial^2 u}{\partial x \partial y} + \frac{1}{2}a_{yy}(\mathbf{x})\frac{\partial^2 u}{\partial y^2} + d_x(\mathbf{x})\frac{\partial u}{\partial x} + d_y(\mathbf{x})\frac{\partial u}{\partial y}}_{\mathcal{G}u} + c(\mathbf{x})u = f(\mathbf{x}) & \text{if } \mathbf{x} \in \Omega, \\ u = g(\mathbf{x}) & \text{if } \mathbf{x} \in \partial\Omega \end{cases}$$
(6)

where $a_{xx}, a_{yy}, a_{xy}, d_x, d_y, c, f, g$ are continuous functions of $\mathbf{x} = (x, y)$, the *reaction coefficient* obeys[4] $c(\mathbf{x}) \leq 0$, and the matrix

$$A(\mathbf{x}) = \begin{bmatrix} a_{xx}(\mathbf{x}) & a_{xy}(\mathbf{x}) \\ a_{xy}(\mathbf{x}) & a_{yy}(\mathbf{x}) \end{bmatrix} = \sigma(\mathbf{x})\sigma^T(\mathbf{x}) \tag{7}$$

is strictly positive definite, so that the *diffusion matrix* $\sigma(\mathbf{x})$ exists.[5] We assume that the coefficients as well as the boundary have the required regularity such that a classical unique solution to (6) exists [19].

Moreover, we will assume that they also enjoy the additional regularity such that the probabilistic representation of the pointwise solution of (6) holds [25]:

$$u(\mathbf{x}_0) = \mathbb{E}[\phi_\tau(\mathbf{x}_0)] := \mathbb{E}\big[ g(\mathbf{X}_\tau)Y_\tau + Z_\tau \,\big|\, \mathbf{X}_0 = \mathbf{x}_0 \big], \tag{8}$$

where the processes $(\mathbf{X}_t, Y_t, Z_t)$ are governed by the following set of SDEs driven by a standard two-dimensional Wiener process $(\mathbf{W}_t)_{t \geq 0}$:

$$\text{while } 0 \leq t \leq \tau(\mathbf{x}_0) \begin{cases} d\mathbf{X}_t = \mathbf{d}(\mathbf{X}_t)dt + \sigma(\mathbf{X}_t)d\mathbf{W}_t, & \mathbf{X}_0 = \mathbf{x}_0 \in \Omega, \\ dY_t = c(\mathbf{X}_t)Y_t dt, & Y_0 = 1, \\ dZ_t = -f(\mathbf{X}_t)Y_t dt + Y_t \mathbf{F}^T(\mathbf{X}_t)d\mathbf{W}_t, & Z_0 = 0. \end{cases} \tag{9}$$

Above, $\mathbf{d} := (d_x, d_y)$ is the drift, $\tau(\mathbf{x}_0) := \inf_{t \geq 0}\{\mathbf{X}_t \in \partial\Omega\}$ is the *first exit time* (or first passage time) of the diffusion starting from a point $\mathbf{x}_0 \in \Omega$ and driven by the differential generator $\mathcal{G}$ defined in (6); which takes place at the *first exit point* $\mathbf{X}_\tau \in \partial\Omega$. Each trajectory yields a *score* $\phi_\tau(\mathbf{x}_0) \in \mathbb{R}$ (the value of the Feynman-Kac functional (8)). Finally, the vector field $\mathbf{F}$ is rather arbitrary, and does not affect the expected value of $\phi_\tau(\mathbf{x}_0)$—and hence of $u(\mathbf{x}_0)$)—but it does affect its variance. (This will be later exploited for variance reduction.)

## 2.3 Numerics of absorbed diffusions

Given $\mathbf{x}_0 \in \Omega$, a numerical solution to (8) consists in a timestepping scheme for (9) with a timestep $0 < h \ll 1$, which yields an approximation to the realisation[6] $\phi_\tau^{(h,j)} \approx \phi_\tau(\omega_j)$ (where $\omega_j$ is the "chance variable" labelling the $j^{th}$ random score), along with an estimator of $\mathbb{E}[\phi_\tau]$, such as the mean. Let $\phi_\tau^{(h)}(\mathbf{x}_0)$ represent the approximate distribution of Feynman-Kac scores constructed by timestepping from $\mathbf{x}_0$. The total pointwise error $e_{N,h}(\mathbf{x}_0)$ is a stochastic variable as well, asymptotically bounded as follows [22]. As $N \to \infty$ and $h \to 0^+$, $e_{N,h}(\mathbf{x}_0)$

---
[4]In spite of this limitation, a class of Helmholtz's equations have been solved in [26].
[5]It can be extracted by performing the Cholesky decomposition of $A$.
[6]We may drop the reference to point $\mathbf{x}_0$ from the random variable $\phi_\tau(\mathbf{x}_0)$ or realisations thereof.

is bounded with Gaussian probabilities $P_q \approx \{0.68, 0.95, 0.997\}$ for $q = \{1, 2, 3\}$ as

$$e_{N,h}(\mathbf{x}_0) := \left| u_{exact}(\mathbf{x}_0) - \frac{1}{N} \sum_{j=1}^{N} \phi_\tau^{(h,j)}(\mathbf{x}_0) \right| \leq \sqrt{\mathbb{E}\left[\left(\phi_\tau^{(h)}(\mathbf{x}_0) - \phi_\tau(\mathbf{x}_0)\right)^2\right]}$$

$$\leq B(\mathbf{x}_0) h^\delta + q \sqrt{\frac{\mathbb{V}\left[\phi_\tau^{(h)}(\mathbf{x}_0)\right]}{N}}, \qquad (10)$$

where $B(\mathbf{x}_0) > 0$ is the *bias constant* (which depends on $\mathbf{x}_0$, the coefficients in (9), the integration domain, and the numerical integrator—but not on $h$); $\delta > 0$ is the *weak convergence rate*, which depends on the integrator, as well as on the smoothness of the SDE coefficients and $\Omega$; and $\mathbb{V}\left[\phi_\tau^{(h)}(\mathbf{x}_0)\right]$ is the variance of the scores, which can be estimated in computing time by the sample variance.

The right hand side of (10) shows the two components of the error: the bias $B(\mathbf{x}_0) h^\delta$, which vanishes as $h \to 0^+$, and the statistical error, which goes to zero as $N \to \infty$. If $\epsilon$ is the pointwise accuracy tolerance, both errors terms must be balanced against one another for the sake of efficiency, for instance at $\epsilon/2$ each. Provided that $B(\mathbf{x}_0)$ is known ($\delta$ is deemed theoretically known), this yields the optimal $h$ for a given error objective. The optimal $N$ follows from letting the statistical error make up the remaining $\epsilon/2$.

The cost of approximating $\mathbb{E}[\phi_\tau(\mathbf{x}_0)]$—measured in evaluations of (9)—is the number of trajectories times the average hops to the boundary, i.e. $N \times (\mathbb{E}[\tau(\mathbf{x}_0)]/h) \propto N/h$. For a balanced simulation with accuracy tolerance $\epsilon$, $N \propto \mathbb{V}[\phi_\tau(\mathbf{x}_0)] \times \epsilon^{-2}$, while $h \propto \epsilon^{1/\delta}$, so that the total complexity is $O(\epsilon^{-2-1/\delta})$. In order to make this cost manageable, high-order SDE weak integrators (for which $\delta = 1$ or perhaps slightly better) as well as statistical error reduction methods are needed. (Both improvements can be achieved simultaneously with Giles' Multilevel Monte Carlo method [20].)

**Naïve Euler-Maruyama scheme, and Gobet-Menozzi shift.** The most basic general-purpose integrator for (9) is explained next. A continuous trajectory $\mathbf{X}_t(\omega_j)$, $(0 \leq t \leq \tau(\omega_j))$, is approximated by the points $\{\mathbf{x}_0 = \mathbf{X}_0^{(h,j)}, \mathbf{X}_1^{(h,j)}, \ldots\}$, drawn according to the Euler-Maruyama scheme [22], until overshooting the boundary at step $K$—so that $\mathbf{X}_K^{(h,j)} \notin \overline{\Omega}$. Then, $\tau(\omega_j) \approx \tau^{h,j}$, where $(K-1)h < \tau^{h,j} < Kh$, and the first exit point for that trajectory is "naïvely" taken as $\mathbf{X}_K^{(h,j)}\big|_{\partial\Omega}$ (the normal projection of $\mathbf{X}_K^{(h,j)}$ on the boundary). The approximate score is $\phi_\tau^{(h,j)} = g\left(\mathbf{X}_K^{(h)}\big|_{\partial\Omega}\right) Y_K^{(h,j)} + Z_K^{(h,j)}$. ($\{Y_k, Z_k\}$ are the obvious counterparts of $\{\mathbf{X}_k\}$.) This scheme is very direct, but has a weak convergence rate of just $\delta \leq 1/2$ [11].

The method introduced by Gobet and Menozzi [21] simply offsets the domain boundary *inwards* by a distance $-0.5826 \|\sigma(\mathbf{X}_k|_{\partial\Omega}) \mathbf{N}(\mathbf{X}_k|_{\partial\Omega})\|_2 h^{1/2}$, where $\mathbf{N}$ is the normal vector to $\partial\Omega$. This scheme is as simple as the naïve one, but enjoys linear weak order $\delta = 1$. It is arguably the best for purely absorbing (i.e. Dirichlet) boundaries [7].

With either scheme (and actually with almost any scheme), it is very difficult to bound $B(\mathbf{x}_0)$ in advance. Multilevel Monte Carlo—which can also incorporate the Gobet-Menozzi shift [20]—is the exception. Unfortunately, it cannot be directly used in PDDSparse, for reasons that will be clear later.

**Variance reduction.** By Itô's isometry, the variance of the Feynman-Kac score can be written as [25]

$$\mathbb{V}\left[\phi_\tau(\mathbf{x}_0)\right] = \mathbb{E}\left[\int_0^{\tau(\mathbf{x}_0)} Y_t^2 \|\sigma^T(\mathbf{X}_t)\nabla u_{exact}(\mathbf{X}_t) + \mathbf{F}(\mathbf{X}_t)\|_2^2 \, dt\right]. \quad (11)$$

Therefore, the variance vanishes if $\mathbf{F} \equiv -\sigma^T \nabla u_{exact}$—one single realisation of the SDE would yield $u_{exact}(\mathbf{x}_0)$. Obviously, this is impossible in practice, since it requires knowing $\nabla u_{exact}$ in all of $\Omega$, and exact integration of (11).

Nonetheless, (11) can be used very effectively in a "multigrid" fashion: first, a rough, but relatively unexpensive, nodal solution $\vec{u}_0$ is obtained with $\mathbf{F} = 0$, leading (after solving the subdomains) to a coarse global approximation $u_0(\mathbf{x})$. The second, more accurate run then takes $\mathbf{F} = -\sigma^T \nabla u_0$ and makes do with fewer trajectories (for a given statistical error target) thanks to the reduced variance. The aggregated wallclock time of these two iterations (or more, see [6]) can be significantly less than one single, accurate solve with $\mathbf{F} = 0$.

## 3 Formulation of PDDSparse for 2D elliptic BVPs

### 3.1 PDDSparse discretisation

We shall now introduce some nomenclature, using Figure 2 as a reference. PDDSparse starts by embedding the actual domain $\Omega \subset \mathbb{R}^2$ in a grid of (open) squares. Such squares are denoted as $\square_k$, with $k = 1, \ldots, M_s$. Note that a square may be contained within $\Omega$, intersect $\partial\Omega$, or lie entirely outside.

To simplify the discussion and the notation, we shall make two assumptions:

**A1.** The boundary $\partial\Omega$ may either i) contain a square side entirely, or ii) intersect it at one single point at most, which moreover is not a grid crossing.

**A2.** The boundary cannot intersect / contain the four sides of any grid square.

The ingredients of the PDDSparse discretisation are:

- *Subdomains.* The $k^{th}$ subdomain is defined as the open, connected set

$$\Omega_k = \Omega \cap \square_k. \quad (12)$$

Thus, PDDSparse subdomains are nonoverlapping by construction.[7]

The BVP restricted to $\Omega_k$ is solved on the last stage of PDDSparse, after the solution of the global BVP on the portion of $\partial\Omega_k$ inside $\Omega$ (called *fictitious interfaces*) has been found. The latter are discretised next.

---

[7]We allow for empty subdomains, if the associated square is fully outside of $\Omega$.

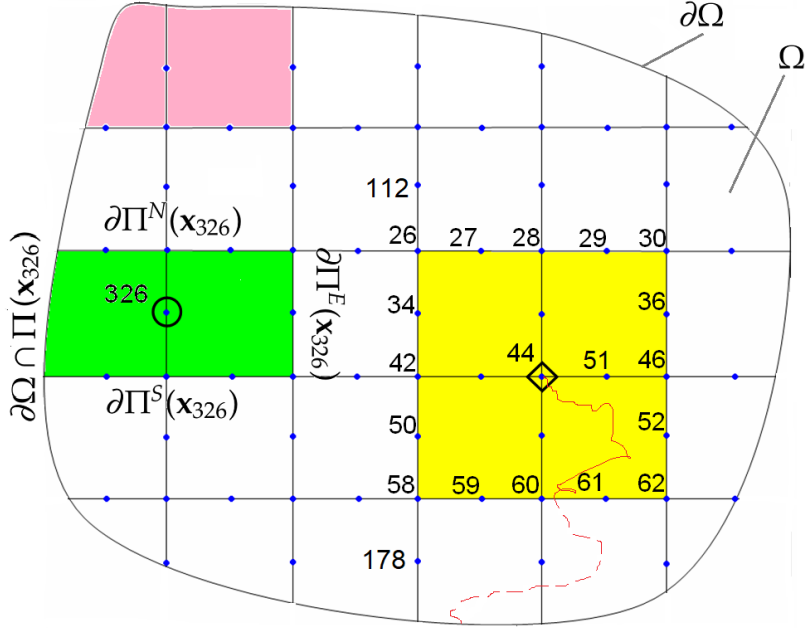

Figure 2: PDDSparse embeds the domain $\Omega$ into a grid of squares from which subdomains and interfaces are defined. Every interfacial knot is associated to one patch; whether a perimeter patch (such as the green one for knot 326, or a floating patch (such as the yellow one for 44 and 51). The stencils of both are: $\mathcal{S}^E_{44} = \mathcal{S}^E_{51} = \{30, 36, 46, 52, 62\}$, $\mathcal{S}^N_{44} = \mathcal{S}^N_{51} = \{26, 27, 28, 29, 30\}$, $\mathcal{S}^W_{44} = \mathcal{S}^W_{51} = \{178, 58, 50, 42, 34, 26, 112\}$, and $\mathcal{S}^S_{44} = \mathcal{S}^S_{51} = \{58, 59, 60, 61, 62\}$. (Note that $\mathcal{S}^E_{44}$ has been extended on both directions.) The red stochastic trajectory exits the patch at $\mathbf{X}_{\tau^{\#}_{44}}$ between $\mathbf{x}_{60}$ and $\mathbf{x}_{61}$; the dashed part is not computed.

- *Knots*, denoted by $\mathbf{x}_i \in \overline{\Omega}$, $i = 1, \ldots, n$, are points on the subdomain's fictitious interfaces where the pointwise solution to the global BVP are calculated by probabilistic methods (or where they take on the Dirichlet BC). The pointwise solution on knot $\mathbf{x}_i$ is the $i^{th}$ *nodal solution*. Knots sit on the lines of the square grid—including at the crossings.

- Every knot $\mathbf{x}_i$ is associated to one *slab* $Q(\mathbf{x}_i)$, defined as the union of all grid squares to which $\mathbf{x}_i$ belongs (actually, to the square boundary):

$$Q(\mathbf{x}_i) = \bigcup_{\mathbf{x}_i \in \partial \Box_k} \overline{\Box}_k. \tag{13}$$

Depending on whether knot $\mathbf{x}_i$ sits at a crossing or between two of them, its slab contains two squares (and is a rectangle), or four of them (and is a square slab). In either case, the four sides of the boundary of $Q(\mathbf{x}_i)$ are

denoted as

$$\partial Q(\mathbf{x}_i) = \partial^E Q(\mathbf{x}_i) \cup \partial^N Q(\mathbf{x}_i) \cup \partial^W Q(\mathbf{x}_i) \cup \partial^S Q(\mathbf{x}_i) =: \bigcup_{O \in \{E,N,W,S\}} \partial^O Q(\mathbf{x}_i), \quad (14)$$

where $O \in \{E, N, W, S\}$ intuitively stand for East, North, West and South.

- The *patch* associated to a knot is the portion of its slab contained in $\overline{\Omega}$:

$$\Pi(\mathbf{x}_i) = Q(\mathbf{x}_i) \cap \overline{\Omega}. \quad (15)$$

The mapping between knots and patches/slabs is not bijective: several knots will belong to the same patch. In order to simplify the notation, we shall write just $Q$ and $\Pi$ when the knot is unimportant for the discussion.

Patches can be *floating*, when they are fully contained in the domain, or *perimeter*, when they intersect the boundary. By assumption **A1**, a floating patch does not intersect $\partial\Omega$. The knot $\mathbf{x}_i$ is called a floating (respectively, perimeter) knot if $\Pi(\mathbf{x}_i)$ is a floating (respectively, perimeter) patch.

- For either kind of patch, let us define the *segments* as

$$\partial^O \Pi = \begin{cases} \partial^O Q(\mathbf{x}_i) & \text{if } \Pi(\mathbf{x}_i) \text{ is a floating patch,} \\ \partial\Omega \cap \partial^O Q(\mathbf{x}_i) & \text{if } \Pi(\mathbf{x}_i) \text{ is a perimeter patch.} \end{cases}$$

By assumption **A1**, a perimeter segment $\partial\Omega \cup \partial^O Q$ is always either an empty set, or a connected set containing more than one single point.

The boundary of a generic patch is then

$$\partial\Pi(\mathbf{x}_i) = \left(\partial\Omega \cup \Pi(\mathbf{x}_i)\right) \cup \left(\bigcup_{O \in \{E,N,W,S\}} \partial^O \Pi(\mathbf{x}_i)\right). \quad (16)$$

The first term above is nonempty only if $\Pi(\mathbf{x}_i)$ is a perimeter patch whose slab overshoots $\partial\Omega$.

- *Stencils.* With each nonempty segment $\partial^O \Pi$ we associate the set of knots sitting on that segment, possibly augmented with knots on either side of the segment ends. Formally, the stencil of segment $\partial^O \Pi$ is the set of knots

$$\mathcal{S}[\partial^O \Pi] = \{\mathbf{x}_j\} \subset \{\mathbf{x}_1, \ldots, \mathbf{x}_n\} \quad (17)$$

(or of knot indices) which verify: *i)* if $\mathbf{x}_i \in \partial^O \Pi$, then $\mathbf{x}_i \in \mathcal{S}[\partial^O \Pi]$ (i.e. the segment knots always belong to the segment stencil). *ii)* All stencil knots lie along the same line. *iii)* If $\mathbf{x}_i \in \mathcal{S}[\partial^O \Pi]$ and $\mathbf{x}_j \in \mathcal{S}[\partial^O \Pi]$; and knot $\mathbf{x}_k$ is sitting in between, then $\mathbf{x}_k \in \mathcal{S}[\partial^O \Pi]$ (i.e. there are no gaps).

The reason why the knots in $\mathcal{S}[\partial^O \Pi]$ must be collinear is to ensure that the solution of the global PDE along them defines a smooth function, and is thus amenable to RBF interpolation.

In practice, we shall talk about *knot stencils*, which are the segment stencils of the given knot's patch. Therefore, let us introduce the notation

$$\mathcal{S}(\mathbf{x}_i) := \mathcal{S}_i := \bigcup_{O \in \{E,N,W,S\}} \mathcal{S}_i^O := \bigcup_{O \in \{E,N,W,S\}} \mathcal{S}[\partial^O \Pi(\mathbf{x}_i)]. \tag{18}$$

## 3.2 A domain decomposition Feynman-Kac formula

Consider an elliptic operator $\mathcal{G}$ as in (6) and the interfacial node $\mathbf{x}_i$. By ellipticity, the nodal solution $u(\mathbf{x}_i)$ is the solution at $\mathbf{x}_i$ of the patch-restricted BVP:

$$\begin{cases} \mathcal{G}u + cu = f & \text{if } \mathbf{x} \in \Pi(\mathbf{x}_i), \\ u = u|_{\partial^E \Pi(\mathbf{x}_i)} & \text{if } \partial^E \Pi(\mathbf{x}_i) \neq \emptyset \text{ and } \mathbf{x} \in \partial^E \Pi(\mathbf{x}_i), \\ u = u|_{\partial^N \Pi(\mathbf{x}_i)} & \text{if } \partial^N \Pi(\mathbf{x}_i) \neq \emptyset \text{ and } \mathbf{x} \in \partial^N \Pi(\mathbf{x}_i), \\ u = u|_{\partial^W \Pi(\mathbf{x}_i)} & \text{if } \partial^W \Pi(\mathbf{x}_i) \neq \emptyset \text{ and } \mathbf{x} \in \partial^W \Pi(\mathbf{x}_i), \\ u = u|_{\partial^S \Pi(\mathbf{x}_i)} & \text{if } \partial^S \Pi(\mathbf{x}_i) \neq \emptyset \text{ and } \mathbf{x} \in \partial^S \Pi(\mathbf{x}_i), \\ u = g & \text{if } \Pi(\mathbf{x}_i) \cap \partial \Omega \neq \emptyset \text{ and } \mathbf{x} \in \Pi(\mathbf{x}_i) \cap \partial \Omega. \end{cases} \tag{19}$$

(Note that the PDE also holds on the fictitious interfaces of the patch.) In order to apply the Feynman-Kac formula, we seek to replace $u|_{\partial^O \Pi(\mathbf{x}_i)}$ by a linear interpolator across the nodal values in that segment's stencil, i.e. $\{u(\mathbf{x}_j) \mid \mathbf{x}_j \in \mathcal{S}_i^O\}$. Let $\tau_i^\#$ be the first-exit time from $\Pi(\mathbf{x}_i)$ of the diffusion $(\mathbf{X}_t)_{t \geq 0}$ governed by the differential generator $\mathcal{G}$, and starting off at $\mathbf{X}_0 = \mathbf{x}_i$. Then, using the compact notation introduced previously, the Feynman-Kac formulas (8)-(9) yield

$$u(\mathbf{x}_i) = \mathbb{E}\left[Z_{\tau_i^\#} \mid \mathbf{X}_0 = \mathbf{x}_i\right] + \mathbb{E}\left[Y_{\tau_i^\#} g\left(\mathbf{X}_{\tau_i^\#}\right) \mathbb{1}_{\{\mathbf{X}_{\tau_i^\#} \in \partial \Omega\}} \mid \mathbf{X}_0 = \mathbf{x}_i\right]$$
$$+ \sum_{O \in \{E,N,W,S\}} \mathbb{1}_{\partial^O \Pi(\mathbf{x}_i) \neq \emptyset} \mathbb{E}\left[Y_{\tau_i^\#} \mathbb{1}_{\{\mathbf{X}_{\tau_i^\#} \in \partial^O \Pi(\mathbf{x}_i)\}} u|_{\partial^O \Pi(\mathbf{x}_i)}(\mathbf{X}_{\tau_i^\#}) \mid \mathbf{X}_0 = \mathbf{x}_i\right]. \tag{20}$$

The above expression is purely formal since the interfacial BCs $u|_{\partial^O \Pi(\mathbf{x}_i)}$ are unknown. Using RBFs[8], it turns out that the "Dirichlet BC" can be written as

$$u|_{\partial^O \Pi(\mathbf{x}_i)}(\mathbf{z}) \approx \sum_{\mathbf{x}_j \in \mathcal{S}_i^O} u(\mathbf{x}_j) H_{ij}^O(\mathbf{z}) \qquad \text{if } \mathbf{z} \in \partial^O \Pi(\mathbf{x}_i), \tag{21}$$

where $H_{ij}^O(\mathbf{z})$ is the $j^{th}$ cardinal function derived in (5), taking the collinear set of knots in $\mathcal{S}_i^O$ as RBF pointset. Then,

$$u(\mathbf{x}_i) \approx \mathbb{E}\left[Z_{\tau_i^\#} \mid \mathbf{X}_0 = \mathbf{x}_i\right] + \mathbb{E}\left[Y_{\tau_i^\#} g\left(\mathbf{X}_{\tau_i^\#}\right) \mathbb{1}_{\{\mathbf{X}_{\tau_i^\#} \in \partial \Omega\}} \mid \mathbf{X}_0 = \mathbf{x}_i\right]$$
$$+ \sum_{O \in \{E,N,W,S\}} \mathbb{1}_{\{\partial^O \Pi(\mathbf{x}_i) \neq \emptyset\}} \sum_{\mathbf{x}_j \in \mathcal{S}_i^O} \mathbb{E}\left[Y_{\tau_i^\#} \mathbb{1}_{\{\mathbf{X}_{\tau_i^\#} \in \partial^O \Pi(\mathbf{x}_i)\}} H_{ij}^O(\mathbf{X}_{\tau_i^\#}) \mid \mathbf{X}_0 = \mathbf{x}_i\right] u(\mathbf{x}_j). \tag{22}$$

---

[8]Any other interpolation scheme would also be possible, as long as it is linear.

Denoting as $u_k \approx u(\mathbf{x}_k)$ $(k = 1, \ldots, n)$ the resulting approximate nodal solution (with interpolation error), and rearranging terms, the system arises:

$$G(u_1, \ldots, u_n)^T = \vec{b}, \quad \text{or} \quad G\vec{u} = \vec{b}, \tag{23}$$

whose entries are

$$G_{ij} = \begin{cases} 1 & \text{if } i = j, \\ -\mathbb{E}\left[Y_{\tau_i^\#}\mathbb{1}_{\{\mathbf{X}_{\tau_i^\#} \in \partial^O \Pi(\mathbf{x}_i)\}} H_{ij}^O(\mathbf{X}_{\tau_i^\#}) \,\big|\, \mathbf{X}_0 = \mathbf{x}_i\right] & \text{if } \mathbf{x}_j \in \mathcal{S}_i^O, \, O \in \{E, N, W, S\}, \\ 0 & \text{otherwise.} \end{cases} \tag{24}$$

and

$$b_i = \mathbb{E}\left[Z_{\tau_i^\#} \,\big|\, \mathbf{X}_0 = \mathbf{x}_i\right] + \mathbb{E}\left[Y_{\tau_i^\#} g\left(\mathbf{X}_{\tau_i^\#}\right) \mathbb{1}_{\{\mathbf{X}_{\tau_i^\#} \in \partial\Omega\}} \,\big|\, \mathbf{X}_0 = \mathbf{x}_i\right]. \tag{25}$$

### 3.3 Discussion

The key advantage of PDDSparse over PDD is keeping the Feynman-Kac diffusions confined within patches. Since this takes only the PDE / SDE coefficients on the given patch, the "memory issue" of PDD is fixed—the subdomain size (hence the patch's) can be adapted to fit the memory of the supercomputer's individual cores / nodes.

PDDSparse paves the way for using memory-limited, unexpensive GPUs—ideal for SDEs—in order to solve BVPs on a truly massive scale. (See Table 1 in this regard.) Except for the final linear system, PDDSparse remains essentially fault-tolerant. Indeed, the loss of a batch of trajectories upon failure of a core is negligible and immediately replaceable.

The required Monte Carlo simulations are orders of magnitude faster / more accurate with PDDSparse than they were with PDD. The cost of constructing the nonzero entries of the $i^{th}$ equation of PDDSparse equals the computational cost of integrating $N_i$ bounded diffusions from knot $\mathbf{x}_i$ with Euler-Maruyama. As discussed in Section 2.3, that cost is proportional to $\mathbb{E}[\tau^\#(\mathbf{x}_i)]$. The expectation of first-exit times from inside $\theta \subset \mathbb{R}^2$ obeys the BVP [8]

$$\mathcal{G}\mathbb{E}[\tau^\theta(\mathbf{x})] + 1 = 0 \quad \text{if } \mathbf{x} \in \theta, \qquad \mathbb{E}[\tau^\theta(\mathbf{x})] = 0 \quad \text{if } \mathbf{x} \in \partial\theta. \tag{26}$$

It can be proved by inspection that, if $\Lambda > 0$ and $\theta' = \{\mathbf{x}' = \Lambda\mathbf{x} \,|\, \mathbf{x} \in \theta\}$, then $\mathbb{E}[\tau^{\theta'}(\Lambda\mathbf{x})] = \Lambda^2 \mathbb{E}[\tau^\theta(\mathbf{x})]$. In words, the mean number of hops to the boundary is proportional to the area. Consequently, if $\Omega$ is a square of area $|\Omega|$ partitioned into $M_{subs} = H^2$ subdomains[9] (of area $|\Omega|/H^2$), the Feynman-Kac trajectories take $\mathcal{O}(M_{subs}^{-1})$ of the wallclock time in PDDSparse than they would with PDD.

Moreover, the statistical error of the Feynman-Kac scores will also typically decrease with the area of the actual domain of integration. To justify why, note that (whether $\mathbf{F} = 0$ in (11) or not), the pointwise variance

$$\mathbb{V}[\phi_\tau(\mathbf{x}_0)] = \mathbb{E}\left[\int_0^{\tau^\theta(\mathbf{x}_0)} Y_t^2 \|\sigma^T \nabla u_{ex}\|_2^2 \, dt \,\bigg|\, \mathbf{X}_0 = \mathbf{x}_0\right] \tag{27}$$

---

[9]We shall stick to this setup of $\Omega$ for the rest of this section, check Figure 1.

is a nondecreasing function of the first exit time. For instance, for a Poisson equation ($Y_t \equiv 1$) whose solution is relatively homogeneous across the patches (such as the one discussed in Section 5), the integral (27) grows roughly linearly with $\mathbb{E}[\tau^\theta(\mathbf{x}_0)]$. Since the number of trajectories (from a knot and to within a given $\epsilon$) is proportional to the variance (at that knot), in our gridded square $\Omega$ example the cost $C$ goes as $C \propto \mathbb{E}[\tau^\Omega(\mathbf{x}_0)] \times \mathbb{E}[\tau^\Omega(\mathbf{x}_0)] \propto M_{subs}^{-2}$, or $C \propto H^{-4}$.

It makes intuitively sense that the bias constant $B(\mathbf{x}_0)$ in (10) increases with the average flight time, since the Euler-Maruyama trajectories will tend to deviate more from their continuous SDE. (For the dependence of the bias constant on the terminal time in the case of unbounded parabolic diffusions, see [3, 4].)

Also with regard to the bias, it is true that PDDSparse trajectories are integrated within rectangles—regardless of the smoothness of $\partial\Omega$—and corners are well known to spoil the accuracy and weak convergence rate of SDE integrators [7]. It turns out that this issue can be easily fixed by enlarging the segment stencils beyond the corners, as shown in Figure 3. If a discrete trajectory overshoots the patch boundary diagonally, the normal projection onto the patch boundary is not defined, and the corner should be taken as first-exit point. However, in our setting, the solution to the original BVP also exists in the region outside the patch, so that projection can be done on the closest extended stencil. This prevents the corner knots from accummulating spurious first-exit points, and becoming oversampled in the RBF interpolant.

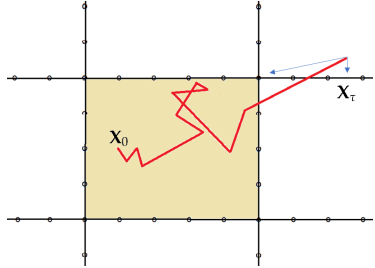

Figure 3: Extending the stencils beyond the patch sides allows for more accurate treatment of SDE "diagonal" boundary overshoots (than projecting on corners).

In sum, both the bias and the statistical error of the SDE simulations are no worse, but very likely significantly better, in PDDSparse than in PDD.

After constructing $G$ and $\vec{b}$ in an embarrassingly parallel way, the system $G\vec{u} = \vec{b}$ must be solved, and this task is no longer fully parallelisable, owing to Amdahl's law. The PDDSparse system, though, is of the size of the interfacial knots only, i.e. $n = O(\sqrt{M})$, where $M$ are the number of point values required to discretise *the whole solution*—the fine grid in Figure 1.

Interpolation errors are negligible, at least as long as $u(\mathbf{x})$ is smooth. This is due to the fact that RBF interpolation of $C^\infty$ functions is spectrally convergent.

In fact, thanks to this property, it makes sense to use a much lower linear density of nodes on the artificial interfaces than inside the subdomains. For the sake of illustration, let $\Omega$ be a square $[0, H]^2$ partitioned into $H^2$ square subdomains of side length 1, each one discretised into a grid of finite differences with $m^2$ nodes (see Figure 1, left). The number of DoFs of the entire problem is $M = H^2 m^2$, and the accuracy typically $O(m^{-2})$. If the same linear density was used on the PDDSparse interfaces, matrix $G$ would be of size $\approx 2Hm$, i.e. $O\left(M^{\frac{1}{2}}\right)$. However, since in PDDSparse the interfacial values are interpolated with RBFs—instead of finite differences—a much lower linear density than $1/m$ would be perfectly adequate, such as $1/m^\nu$ with $\nu = 2, 3, \ldots$. In that case, matrix $G$ would be only of size $O\left(m^{\frac{1}{2\nu}}\right)$, which is better than by a straightforward algebraic substructuring (Schur complement-like) scheme.

Since floating subdomains are square, subdomain-restricted BVPs can be solved with the pseudospectral method [31] (as it is the case in our implementation). Perimeter subdomains are tackled, during the construction of the system, by Euler-Maruyama-like integrators that only rely on knowledge of the distance to the boundary. For practical purposes, PDDSparse is thus a meshless method, free from the costly generation of well-behaved finite element (FE) meshes and subdomain partitions that are a necessary, and often very time consuming, preprocessing step of large-scale FE simulations.

The last, but potentially important, source of error is related to the stability of the system $G\vec{u} = \vec{b}$. Our numerical observations so far (and theoretical analysis to be published elsewhere [9]) strongly suggest that the condition number of $G$ grows roughly linearly with the number of interfacial nodes. In our experiments, it was at most $O(100)$—thus effectively negligible with respect to Monte Carlo errors. In addition to it, a simple preconditioner has been constructed in [9], which brings the condition number down to $O(1)$ for the same experiments. (Admittedly, those experiments are still far from the supercomputing regime PDDSparse has been designed to take advantage of.)

## 4 Implementation

### 4.1 Error control

**Heuristic interpretation of nodal errors.** The entries in the $i^{th}$ row of $G$ and $\vec{b}$ are constructed based on $N_i$ stochastic trajectories from knot $\mathbf{x}_i$ with a timestep $h_i$. Then, a PDDSparse simulation is one realisation of the perturbed system

$$(G + \Delta G')(\vec{u}_{ex} + \Delta \vec{u}') = \vec{b} + \Delta \vec{b}', \qquad (28)$$

where $\vec{u}_{ex}$ are the exact nodal solutions, and the perturbations $\Delta G'$ and $\Delta \vec{b}'$ carry bias and statistical error, as well as interpolation error. We are interested in the propagated error $\Delta \vec{u}'$.

Based on Section 3.3, we neglect RBF interpolation errors and assume that the system is well conditioned, or it can be successfully preconditioned. This is a drastic simplification, but has worked well so far and is easy to implement and interpret. Expanding (28), it holds

$$G\Delta\vec{u}' + \Delta G\vec{u} = \Delta\vec{b}' - \Delta G'\Delta\vec{u}'. \tag{29}$$

After computing $(\vec{u}_{ex} + \Delta\vec{u}')$, the solutions at knots $i = 1, \ldots, n$ can be approximated again by integrating a set of Feynman-Kac trajectories inside their patches. This yields, by the very construction of PDDSparse

$$\vec{u}_{ex} + \Delta\vec{u}'' = (\vec{b} + \Delta\vec{b}') + (I - G - \Delta G')(\vec{u}_{ex} + \Delta\vec{u}''), \tag{30}$$

where $(\Delta\vec{u}'')_i$ is now the error of the approximation of a bounded SDE upon integration with the Gobet-Menozzi scheme, timestep $h_i$, and $N_i$ trajectories—whose error structure (10) is well understood. Then, plugging (28) into (29),

$$\Delta\vec{u}'' \approx \Delta\vec{u}'. \tag{31}$$

The point of (31) is that it allows us to estimate the $\{N_i\}$ and $\{h_i\}$ for a given accuracy, by postprocessing a warming-up PDDSparse simulation. Furthermore, it accounts for the success of the variance reduction scheme based on pathwise control variates that will be introduced later in this section.

**Warm-up PDDSparse simulation.** A relatively coarse PDDSparse simulation is first carried out with $N_i = N_0, h_i = h_0$ on all knots, yielding a nodal solution $\vec{u}_0$ and then a global solution (after solving the subdomain-restricted BVPs) $u_0(\mathbf{x})$.

The coarse global solution $u_0$ allows for an estimation of pointwise bias constants and variances according to heuristic (31), as explained next.

**Estimation of nodal biases.** A posteriori estimates of the bias coefficient of a Monte Carlo simulation are notoriously expensive to produce, typically calling for at least two simulations with different timesteps. Here, we present an idea which takes one single set of simulations at timestep $h_0'$. Let $\phi_\tau^{EM}$ and $\phi_\tau^{GM}$ be respectively the naïve Euler-Maruyama and Gobet-Menozzi approximations to the Feynman-Kac score. By the triangle inequality,

$$\left|\mathbb{E}[\phi_\tau^{EM} - \phi_\tau]\right| = \left|\mathbb{E}[\phi_\tau^{EM} - \phi_\tau^{GM}] + \mathbb{E}[\phi_\tau^{GM} - \phi_\tau]\right|$$
$$\leq \left|\mathbb{E}[\phi_\tau^{EM} - \phi_\tau^{GM}]\right| + \left|\mathbb{E}[\phi_\tau^{GM} - \phi_\tau]\right|. \tag{32}$$

For $h_0'$ small enough, one can assume that

$$\left|\mathbb{E}[\phi_\tau^{GM} - \phi_\tau]\right| < \frac{1}{2}\left|\mathbb{E}[\phi_\tau^{EM} - \phi_\tau]\right|, \tag{33}$$

so that (32) implies $\left|\mathbb{E}[\phi_\tau^{GM} - \phi_\tau]\right| \leq \left|\mathbb{E}[\phi_\tau^{EM} - \phi_\tau^{GM}]\right|$. Moreover, inserting (33) into the "reverse" triangle inequality gives,

$$\left|\mathbb{E}[\phi_\tau^{EM} - \phi_\tau]\right| = \left|\mathbb{E}[\phi_\tau^{EM} - \phi_\tau^{GM}] - \mathbb{E}[\phi_\tau - \phi_\tau^{GM}]\right|$$
$$\geq \left|\left|\mathbb{E}[\phi_\tau^{EM} - \phi_\tau^{GM}]\right| - \left|\mathbb{E}[\phi_\tau^{GM} - \phi_\tau]\right|\right|$$
$$= \left|\mathbb{E}[\phi_\tau^{EM} - \phi_\tau^{GM}]\right| - \left|\mathbb{E}[\phi_\tau^{GM} - \phi_\tau]\right|. \tag{34}$$

Combining (32) and (34) suggests that $\left|\mathbb{E}[\phi_\tau^{EM} - \phi_\tau]\right| \approx \left|\mathbb{E}[\phi_\tau^{EM} - \phi_\tau^{GM}]\right|$, whence

$$\left|\mathbb{E}[\phi_\tau^{GM} - \phi_\tau]\right| \lesssim \frac{1}{2} \left|\mathbb{E}[\phi_\tau^{EM} - \phi_\tau^{GM}]\right| \tag{35}$$

is our final upper bound estimate of the Gobet-Menozzi bias in terms of a computable quantity.

In order to produce a low variance estimate of the previous formula, after the warm-up simulation we draw the two sets of scores $\{\phi_\tau^{GM,j,h'_0}\}_{j=1}^{N'_0}$ and $\{\phi_\tau^{EM,j,h'_0}\}_{j=1}^{N'_0}$ from *one single set of $N'_0$ trajectories*. They are integrated from $\mathbf{x}_i$ within its patch $\Pi(\mathbf{x}_i)$, using $u_0|_{\partial \Pi(\mathbf{x}_i)}$ as Dirichlet BC and a timestep $h'_0$. Every individual trajectory will yields either both scores simultaneously or the GM one first (if it overshoots the shrunken boundary, but not the patch one). In this case, the trajectory will simply be resumed until overshooting the patch boundary, yielding the EM score. Since both sets of scores correspond to totally or partially identical trajectories, we expect for the statistical errors of both means to be highly positively correlated. Consequently,

$$\frac{1}{N'_0} \sum_{j=1}^{N'_0} \left( \phi_\tau^{EM,h'_0,j} - \phi_\tau^{GM,h'_0,j} \right) \approx \mathbb{E}[\phi_\tau^{EM} - \phi_\tau^{GM}] \tag{36}$$

is expected to lead to a low variance estimate of (35) even if $N'_0$ is not too large.

**Estimation of nodal variances.** According to heuristic (31), the sample variances of $\phi_\tau^{GM}(\mathbf{x}_i)$ can be computed using the same trajectories as were used to estimate the biases, for every knot within its patch.

## 4.2 Variance reduction and setting $\{N_i, h_i\}$

In PDDSparse, the statistical error can be substantially reduced at negligible cost thanks to variance reduction based on pathwise control variates. With the coarse global approximation $u_0$, the variable $\mathbf{F} := \sigma^T \nabla u_0$ is constructed and stored on a patch basis, as an interpolatory look-up table. Later, the Itô integral

$$\xi := -\int_0^{\tau_i^\#} Y_t \mathbf{F}(\mathbf{X}_t) \cdot d\mathbf{W}_t \tag{37}$$

can be drawn alongside the scores $\phi_\tau$. In fact, $\xi$ is a control variate of $\phi_\tau$, since

$$\mathbb{E}[\phi_\tau + \xi^h] = \mathbb{E}[\phi_\tau] \qquad \text{but} \qquad \mathbb{V}[\phi_\tau^h + \xi^h] \ll \mathbb{V}[\phi_\tau^h]. \tag{38}$$

The drop in variance increases as $u_0(\mathbf{x}) \to u_{ex}(\mathbf{x})$ and as $h \to 0$. Moreover, nested control variates are possible and efficient within PDD—see [6] for details. The quantity of interest is the variance reduction factor, and it can be proved that

$$\mathbb{V}[\phi_\tau^h + \xi^h] = \mathbb{V}[\phi_\tau]\left(1 - \rho^2(\phi_\tau^h, \xi^h)\right), \tag{39}$$

where $-1 \leq \rho \leq 1$ is Pearson's correlation coefficient between $\phi_\tau^h$ and $\xi^h$.

The idea, then, is to use the same trajectories with $N_0'$ and $h_0'$ (that were used to estimate biases and variances) to also bound Pearson's correlation on a knot by knot basis. It was proved in [6] that $1 - \rho^2(\phi_\tau^h, \xi^h)$ grows linearly with $h$. Therefore, since $h_i \leq h_0 \leq h_0'$, $\rho^2(\phi_\tau^{h_0'})$ leads to a lower bound of the final variance reduction factor.

**Setting the $\{h_i\}_{i=1}^n$ and the $\{N_i\}_{i=1}^n$.** Given an accuracy objective $\epsilon$ we determine the $h_i$ and $N_i$ for each knot according to the error structure (10), splitting the error into two thirds for the statistical component, plus one third for the bias.

For the timestep $h_i$, we employ the bound (35) as a conservative estimate, and let (adding the safeguard that $h_i \leq h_0$, and taking the mean over $N_0'$ below)

$$h_i = \min\left\{h_0, \frac{2h_0'\epsilon}{3\left|\mathbb{E}\left[\phi_\tau^{GM,h_0'}(\mathbf{x}_i) - \phi_\tau^{EM,h_0'}(\mathbf{x}_i)\right]\right|}\right\}. \qquad (40)$$

In which concerns $N_i$, we directly incorporate our forecast of variance reduction (39), and extract it from (with a 95% confidence interval)

$$\frac{2}{3}\epsilon = 2\sqrt{\frac{\mathbb{V}\left[\phi_\tau^{GM,h_0'}\left(1 - \rho^2(\phi_\tau^{GM,h_0'}, \xi^{h_0'})\right)\right]}{N_i}}. \qquad (41)$$

In order to guarantee a minimal number of trajectories, and for load balancing purposes, the $N_i$ from (41) is rounded to the closest value of $k \times N_{job}$, $k \in \mathbb{N}$.

Errors inside subdomains are smaller than nodal errors thanks to ellipticity.

### 4.3 Load balancing and pseudocode

One of the cores is chosen to be the server with the rest of them being workers. In order to optimise the process of constructing the PDDSparse system—which is by far the most computationally expensive part—the $N = \sum_{i=1}^n N_i$ trajectories that have to be computed to fill each row of $G\vec{u} = \vec{b}$ are split in $N/N_{job}$ jobs so that many cores can take part in their computation. The server determines how many jobs have to be done and sends them to the workers, which integrate the trajectories of each job. After finishing a job, the worker sends back to the server the fresh contributions to to $G$ and $\vec{b}$, and perhaps data related to the bias and variance estimation. The server reduces this contributions and solves the linear system. Later, the server sends the interfacial solution to the workers where the solution inside the subdomains is computed in parallel using a pseudospectral method. Judiciously setting $N_{job}$—perhaps dynamically—ensures that no workers are idle during the Monte Carlo phases.

Beyond the Monte Carlo part, the code is not yet optimised. At this moment, the linear system is solved on the server only, and the last-stage, well-posed, patch-restricted BVPs are solved by just one worker each. Clearly many resources are left iddle by the incomplete current code.

In the very near future, both tasks will be parallelised. (Even though they involve linear algebra, and so are subject to Amdahl's law, great speedup can be attained before saturation.) In the simulations presented in Section 5, this would have a modest effect on the overall wallclock time (given that the discretisation of the BVP is not huge), but would very significantly improve strong scalability.

The complete pseudocode is listed below.

---

**Input:** $N_0, h_0, N_{job}$, BVP, domain discretisation, accuracy target $\epsilon$.

**Phase I (Warm up)**

Build $G_0, \vec{b}_0$ with $N_0, h_0$; solve $G_0 \vec{u}_0 = \vec{b}_0$; solve subdomains for $u_0(\mathbf{x})$

**Phase II (Calibration)** Estimate constants (letting $N_0' = N_0$ and $h_0' = h_0$)

Run $N_0$ "mixed GM/EM" trajectories per knot with timestep $h_0$

Bound bias constants according to (35) and (36)

Draw $\xi^{h_0}$ alongside and compute $\rho^2(\phi_\tau^{h_0}, \xi^{h_0})$ for each knot

Estimate nodal sample variances $\mathbb{V}[\phi_\tau] \approx \mathbb{V}[\phi_\tau^{h_0}]$

Set nodal timesteps $\{h_i\}_{i=1}^n$ according to (40)

Set $\{N_i\}_{i=1}^n$ according to (41)

**Phase III (Production)** Final PDDSparse run to precision $\epsilon$

Build $G, \vec{b}$ with $\{N_i, h_i\}_{i=1}^n$; solve $G\vec{u} = \vec{b}$; solve subdomains for $u(\mathbf{x})$

---

## 5 Numerical experiment

### 5.1 Setup

PDDSparse was tested on Galileo100, a supercomputer at the CINECA facility in Italy. Galileo100 is equipped with 528 computational nodes consisting of two Intel CascadeLake 8260 CPUs, with 24 cores each running at 2.4 GHz and 384GB RAM. Nodes are connected through an Infiniband network capable of a maximum bandwidth of 100 Gbit/s. The SLURM partition assigned to our project allowed us to use up to 32 nodes (1536 cores), during a maximum wallclock time of 24 hours per job, and a total allowance of 82000 hours for the project overall.

The algorithm was written in C++/Octave and parallelised with MPI. Octave is only needed when computing the solution inside the subdomains, by the pseudospectral method [31]. The `Eigen3` library was used to manage linear algebra, and `gsl`, to generate random numbers and interpolate look-up tables.

The BVP to be solved is a Poisson equation with Dirichlet BCs

$$\nabla^2 u = f(x, y) \text{ if } \mathbf{x} \in \Omega, \qquad u = u_{ex}(x, y) \text{ if } \mathbf{x} \in \partial\Omega \qquad (42)$$

where $f(x, y)$ is defined such that the solution of the problem is

$$u_{ex} = 3 + \tfrac{1}{3} \sin\left(\sqrt{1 + \tfrac{x^2}{100} + \tfrac{y^2}{50}}\right) + \tfrac{1}{3} \tanh\left[\sin\left(\tfrac{3x}{25} + \tfrac{y}{20}\right) + \sin\left(\tfrac{x}{20} - \tfrac{3y}{25}\right)\right]. \quad (43)$$

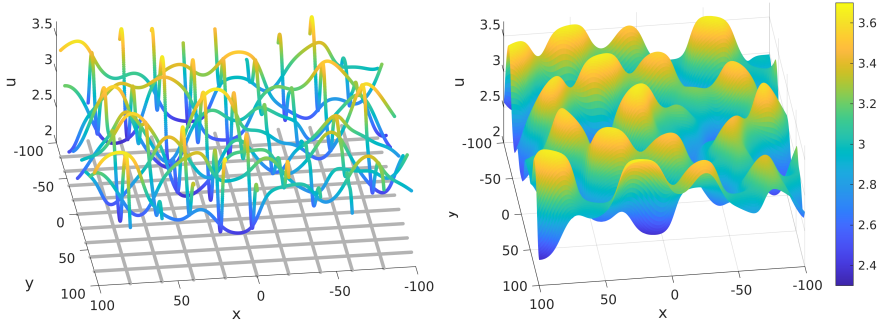

Figure 4: Solution (43) on the complete domain (right) and on the interfacial knots (left), overlaid onto the $10 \times 10$ grid of subdomains.

The domain is the square $\Omega = [-100, 100]^2$, decomposed into $10 \times 10$ subdomains with a total of 409600 DoFs.[10] Each interface is supported by 64 knots. Therefore, there are both rectangular and square patches. Patch side stencils were extended by half a segment on each direction whenever possible.[11] The shape parameter $c$ of the RBF (3) was chosen such that the spectral condition number of the RBF interpolation matrix (2) was $\mathcal{O}(10^{10})$.[12] The discretised interfaces consist of $n = 11277$ knots. The PDDSparse system thus attains a shrinkage of 97.25%, with a 3.35% of nonzero entries[13] and a condition number of ca. 192.[14]

In order to balance the computational load, the number of trajectories per job was set to be $N_{job} = 200$. In the presented simulations, Phases I and II of the algorithm have used $h_0 = h'_0 = 0.08$ and $N_0 = N'_0 = 1000$.

---

[10] Larger problems were also successfully explored. This size was dictated by the need to construct a set of scalability and error curves within the supercomputer access constraints outlined before.

[11] For instance, as it occurs with $\mathcal{S}^W_{44}$ in Figure 2, where knots 112 and 178 have been included in that side's stencil.

[12] Results proved not to be overly sensitive to the shape parameter, though.

[13] While that is too many for $G$ in this example to qualify as an actual sparse matrix, note that the bandwidth is independent of $M_{subs}$, and that PDDSparse is designed for far more subdomains than $10 \times 10$ as here. In a large-scale simulation, sparsity would be orders of magnitude smaller.

[14] We chose not to precondition the PDDSparse matrix in this work, but did check that it can be successfully done [9].

## 5.2 Results

### 5.2.1 Error analysis

We start by demonstrating how the numerical error can be controlled from Phases I through III of PDDSparse, as explained in Section 4. The error incurred in Phase I (warming up) is dominated by the noise, as shown in Figure 5.

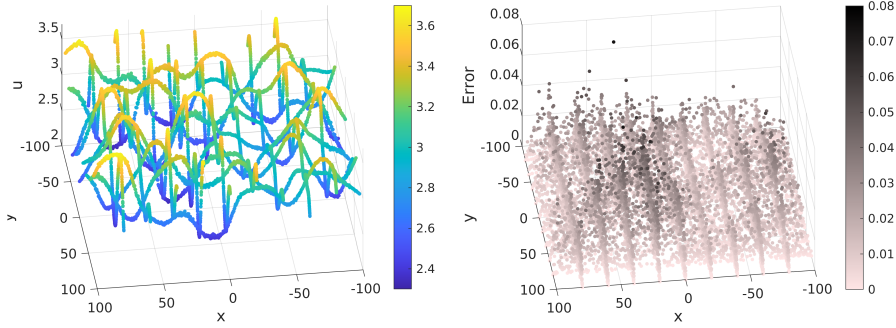

Figure 5: Interfacial solution of Phase I of PDDSparse (left), and its absolute error with respect to (43) (right).

This nodal solution $\vec{u}_0$ leads to the lookup tables of $u_0(x, y)$ that will be used in Phase II (calibration), in order to estimate the $\{N_i, h_i\}$ parameters required to bound the overall numerical error in $\Omega$ by $\epsilon = 0.01$ with ca. 95% probability. In this example, the warmup solution has a maximum error of 0.074 and a RMSE (root mean squared error) of 0.016 on the interfacial knots, with just 48.98% of the knots' solutions having an absolute error lower than $\epsilon$.

In the Phase II (calibration), the variance and bias on all the knots are estimated by taking $u_0(x, y)$ from the first stage as an approximation to the actual solution of the problem. The computation of Pearson's correlation coefficient $\rho\left(\phi_\tau^{h_0 \tau}, \xi^{h_0}\right)$ allows to predict what is going to be the performance of the variance reduction on $\phi_\tau^{h_i}$ during the upcoming Phase III (see Figure 6).

Finally, Phase III (production) is run next. Using the biases and variances estimated in the previous stage, $G$, $\vec{b}$, and $\vec{u}$ are computed again using the control variates and the adjusted values of $N$ and $h$. The maximum error among the interfacial values turns out to be 0.0224 with a RMSE of 0.0044 and a 97.03% of the knots' absolute errors covered by $\epsilon = 0.01$ (Figure 7). The nodal error histograms displayed in Figure 8 attest to the denoising and debiasing achieved between Phases I and III, as well as to the fact—consistent with the heuristic (31)—that nodal errors are reasonably normal (the kurtosis of a Gaussian is 3).

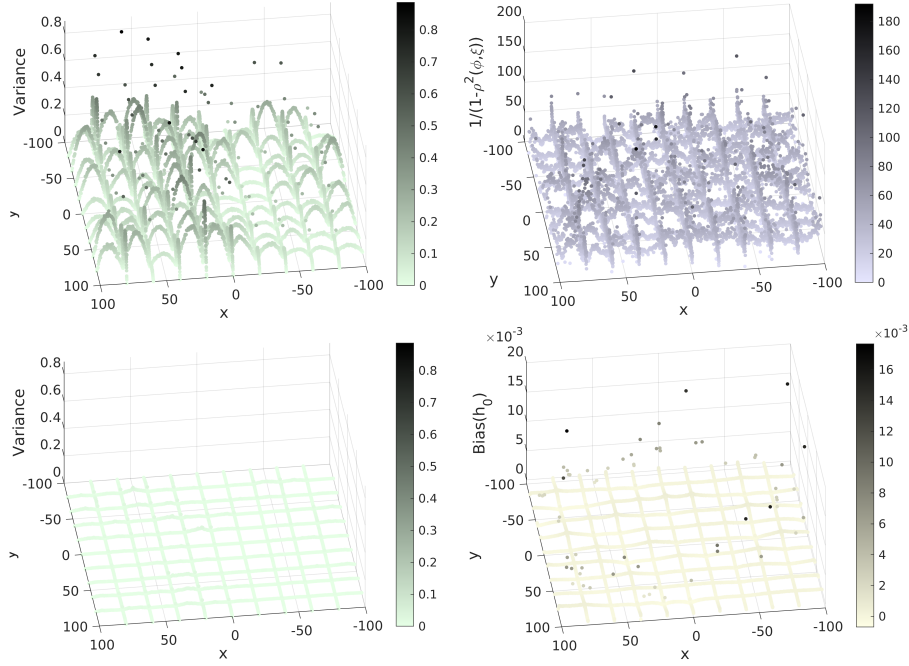

Figure 6: Calibration during Phase II of PDDSparse. The estimated variance of the scores of interfacial knots after Phase I (top left) will be theoretically decreased by a factor of at least $1 - \rho^2\left(\phi_\tau^{h_0}, \xi^{h_0}\right)$ (top right) using pathwise control variates. The realised variance (bottom left) is plotted with the same color scale as the original variance in order to highlight its drop in magnitude. (The average variance reduction factor is 31.37.) On the bottom right, the estimated bias at $h_0$ (it averages $3 \times 10^{-5}$). On most knots there will be no need to take $h_i < h_0$ in Phase III.

### 5.2.2 Strong scalability

In order to study the strong scalability of PDDSparse, we execute the algorithm with identical parameters using $P = \{48, 96, 192, 384, 768, 1536\}$ cores of Galileo100 and measure the wallclock time spent doing so. The resulting strong scalability curves are displayed in Figure (9) for three accuracy targets, as well as the wallclock times segregated by algorithmic task.

It must be stressed that the non-Monte Carlo tasks are not yet optimised. In particular, the PDDSparse systems (in Phase I and III) are for the time being solved on *one single core*—which is why its wallclock time does not fall at all with more cores. While not fully like Monte Carlo simulations—owing to Amdahl's law—linear systems can be parallelised, too—as every state of the art domain decomposition algorithm does. Doing the same in PDDSparse, strong

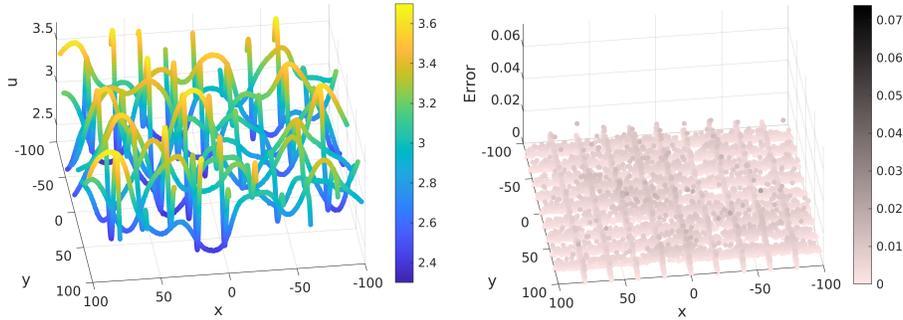

Figure 7: Interfacial solution after Phase III of PDDSparse (left) and its absolute error (right). In order to highlight the drop of error achieved by the intermediate (Phase II) stage, the colormap is the same as that in Figure 5.

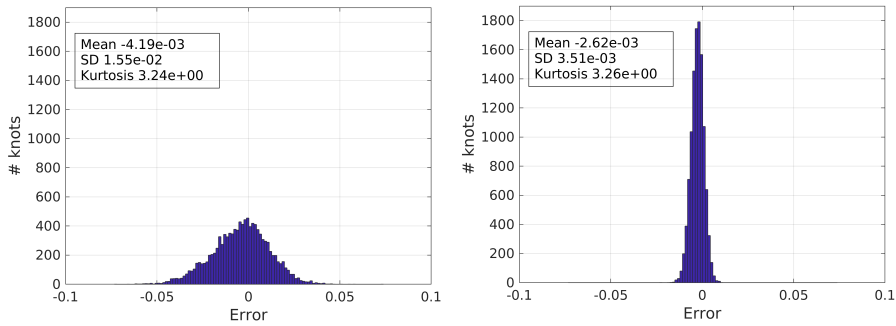

Figure 8: Comparison of the histograms of the error (with sign) of the solution on the interfacial knots after Phase I (left) and after Phase III (right).

scalability would be sustained for far longer than ca. 1500 cores.[15]

Also with the current (provisory) code, each well-posed subdomain is solved by one core only. This means that only 100 cores (as many as subdomains in this experiment) are not idling during that stage. This too can be fixed, leading to slightly shorter wallclock times and better strong scalability.

Finally, Figure 10 plots the empirical acceleration as a function of $P$ and $\epsilon$.

## 5.3 Final remarks

PDDSparse is a hybrid between the state of the art (i.e. algebra based), and probabilistic domain decomposition. The motivation is to recruit more supercomputer cores for the solution of very large BVPs than it is currently possible, by pushing up the limit to strong scalability decreed by Amdahl's curse. The strategy is to insert stochastic calculus and Monte Carlo simulations, which have a better scope for parallelisation than linear algebra.

---

[15]Anyway, that would exceed our access constraints to Galileo100.

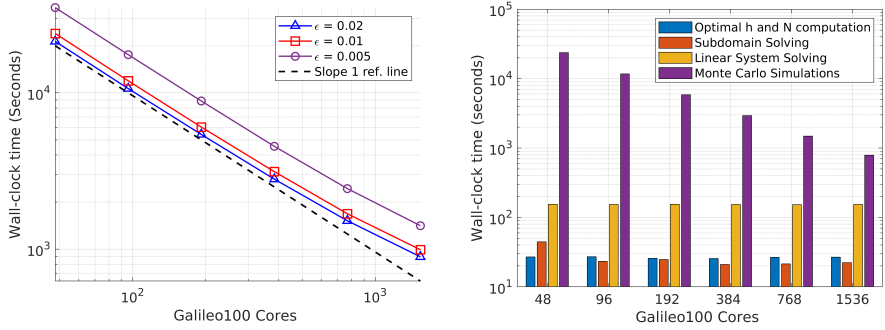

Figure 9: Strong scalability curve of PDDSparse for several accuracy targets $\epsilon$ (left). Right: distribution of times spent on the main tasks of the algorithm.

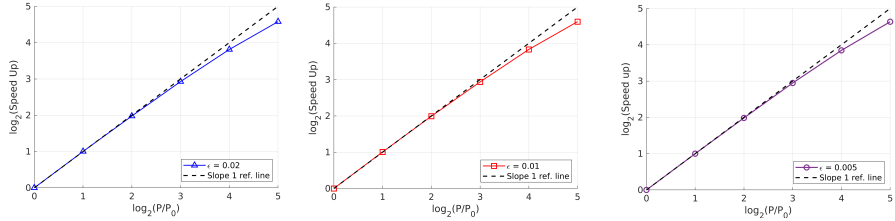

Figure 10: Acceleration curves of PDDSparse with three different values of $\epsilon$. The reference (origin) is the simulation with $P_0 = 48$ cores.

The most obvious advantage of PDDSparse is to shrink the algebraic system down to the interfacial values only, in a direct, sparse, rather general and transparent way. Other features include robustness against core crash and suitability for GPUs. Because the formulation of the interfacial system is so different from the state of the art, it offers new possibilities which deserve further study—including a combination with standard domain decomposition. Extension to three-dimensional elliptic problems is (conceptually) straightforward, as it is to more general boundary conditions [8]. Semilinear elliptic BVPs have been successfully tackled in [26].

Current work is focused on the theoretical investigation on the properties of the PDDSparse system; and on borrowing and adapting existing strategies for solving it in parallel when $n$ is really large—say $n = \mathcal{O}(10^9)$, or $M = \mathcal{O}(10^{18})$ DoFs. Preliminary considerations suggest that the PDDSparse discretisation fits in well with a two-level—or even multigrid—approach. Randomised numerical algebra [24]—increasingly popular among BVP practitioners—is another, fitting avenue.

We have presented proof of concept tests carried out with a prototype code on the Galileo100 supercomputer in CINECA. While strong scalability and error control are satisfactory, wallclock times are still suboptimal. A direct way to improve on them is integrating the SDEs on GPUs. The most current PDDSparse code has a CUDA module for this, which significantly shortens the Monte Carlo timings. For the sake of illustration, consider Table 1, which compares the Monte Carlo time of constructing one row of $G$ and $\vec{b}$ on $\Omega = [-10, 10]^2$, using either a domestic GPU or a CPU.[16]

Table 1: CPU vs. GPU wallclock times, in seconds, with $N = 65536$ trajectories.

| subdomains | CPU | GPU |
|---|---|---|
| 2 x 2 | 1419 | 10.04 |
| 4 x 4 | 359 | 2.57 |
| 8 x 8 | 93 | 0.67 |

# Acknowledgements

FB and JMV were funded by grant 2018-T1/TIC-10914 of the Madrid Regional Government in Spain. JAA thanks the Ministerio de Universidades and specifically the requalification program of the Spanish University System 2021-2023 at the Carlos III University. Access to Galileo100 was awarded thanks to grant HPC173XD3A of the HPC-Europa3 programme of the European Union. A substantial part of this work was carried out during a research stay at the Department of Mathematics and Physics of the Roma Tre University, partially funded by the HPC-Europa3 grant (FB) and a PhD mobility fellowship by Carlos III University of Madrid (JMV). FB and JMV thank Renato Spigler and the rest of the host department for their hospitality, and CINECA and HPC-Europa3 staff for their technical and administrative help. Andrés Berridi is acknowledged for the CUDA code and the GPU vs. CPU comparison.

# References

[1] J. A. Acebrón, M. P. Busico, P. Lanucara, and R. Spigler (2005). *Domain decomposition solution of elliptic boundary-value problems via Monte Carlo and quasi-Monte Carlo methods*. SIAM J. Sci. Comput. 27(2):440–457

[2] S. Badia, A.F. Martín, J. Príncipe (2014) *A highly scalable parallel implementation of balancing domain decomposition by constraints*. SIAM J. Sci. Comput. 36(2):C190–C218

---

[16] The comparison was between a GPU Nvidia Geforce 3080 ti with 10240 cores at 1.67 GHz which cost 1.487.52 €, and a CPU AMD Ryzen 7 3800X with 8 cores at 3.9 GHz and price 247.85 €.

[3] V. Bally and D. Talay (1996) *The Law of the Euler Scheme for Stochastic Differential Equations: I. Convergence Rate of the Distribution Function.* Probab. Th. Rel. Fields 104:43–60

[4] V. Bally and D. Talay (1996) *The Law of the Euler Scheme for Stochastic Differential Equations: II. Convergence Rate of the Density*, Monte Carlo Methods Appl. 2(2):93-128

[5] P. E. Bjørstad and O. B. Widlund (1986) *Iterative methods for the solution of elliptic problems on regions partitioned into substructures.* SIAM J. Numer. Anal., 23(6):1093–1120

[6] F. Bernal and J. A. Acebrón (2016). *A multigrid-like algorithm for probabilistic domain decomposition.* Comput. Math. Appl. 72(7):1790–1810

[7] F. Bernal and J.A. Acebrón (2016) *A comparison of higher-order weak numerical schemes for stochastic differential equations in bounded domains.* Comm. Comput. Phys. 20(3):703–732

[8] F. Bernal (2019) *An implementation of Milstein's method for general bounded diffusions.* J. Sci. Comput. 79(2):867–890

[9] F. Bernal and J. Morón-Vidal (in preparation). *Analysis of a probabilistic domain decomposition algorithm for scientific supercomputing.*

[10] J. F. Bourgat, R. Glowinski, P. Le Tallec, and M. Vidrascu (1988) *Variational formulation and algorithm for trace operator in domain decomposition calculations.* In: T. Chan, R. Glowinski, J. Périaux, and O. Widlund, eds., Domain Decomposition Methods. Second International Symposium on Domain Decomposition Methods, p3–16, Philadelphia, PA, SIAM, 1989

[11] C. Constantini, B. Pacchiarotti and F. Sartoretto (1998) *Numerical approximation for functionals of reflecting diffusion processes.* SIAM J. Appl. Math. 58:73-102

[12] C. R. Dohrmann (2003) *A preconditioner for substructuring based on constrained energy minimization.* SIAM J. Sci. Comput. 25(1):246–258

[13] V. Dolean, P. Jolivet and F. Nataf. An introduction to domain decomposition methods: algorithms, theory, and parallel implementation. Society for Industrial and Applied Mathematics, 2015.

[14] V. Dolean, F. Nataf, R. Scheichl and N. Spillane (2012) *Analysis of a Two-level Schwarz Method with Coarse Spaces Based on Local Dirichlet-to-Neumann Maps.* Computational Methods in Applied Mathematics, 12(4):391-414

[15] M. Dryja (1982) *A capacitance matrix method for Dirichlet problem on polygon region.* Numer. Math., 39:51–64

[16] M. Dryja and O. B. Widlund (1985) *Schwarz methods of Neumann-Neumann type for three-dimensional elliptic finite element problems*. Comm. Pure Appl. Math., 48(2):121–155

[17] C. Farhat, K. Pierson, and M. Lesoinne (2000) *The second generation FETI methods and their application to the parallel solution of large-scale linear and geometrically non-linear structural analysis problems*. Computer methods in applied mechanics and engineering, 184(2-4):333-374

[18] M. Freidlin. Functional integration and Partial Differential Equations. Annals of Mathematics Studies, Vol. 109, Princeton University Press, 1985

[19] D. Gilbarg and N. S. Trudinger. Elliptic partial differential equations of second order. Springer-Verlag, Berlin, New York, 1977

[20] M.B. Giles and F. Bernal (2018) *Multilevel estimation of expected exit times and other functionals of stopped diffusions*. SIAM/ASA Journal on Uncertainty Quantification 6 (4):1454-1474

[21] E. Gobet and S. Menozzi (2010) *Stopped diffusion processes: overshoots and boundary correction*. Stoch. Proc. Appl. 120:130-162

[22] P.E. Kloeden and E. Platen. Numerical Solution of Stochastic Differential Equations. Springer, Applications of Mathematics 23, 1999

[23] J. Mandel (1993) *Balancing domain decomposition*. Comm. Numer. Meth. Engrg., 9: 233–241

[24] P.G. Martinsson, and J.A. Tropp (2020) *Randomized numerical linear algebra: Foundations and algorithms*. Acta Numerica 29:403 - 572

[25] G.N. Milstein and M.V. Tretyakov. Stochastic Numerics for Mathematical Physics. Springer, 2004.

[26] J. Morón-Vidal, F. Bernal, and R. Spigler (submitted) *Iterative schemes for probabilistic domain decomposition*

[27] W. Proskurowski and O. Widlund (1976) *On the numerical solution of Helmholtz's equation by the capacitance matrix method*. Math. of Comp., 30(135):433–468

[28] Y. Saad. Iterative methods for sparse linear systems, 2nd ed.. SIAM, 2003

[29] D. Smith, P. Bjørstad and W. Gropp. Domain Decomposition: Parallel Multillevel Methods for Elliptic Partial Differential Equations. Cambridge University Press, Cambridge (2004)

[30] A. Toselli and O. Widlund, Domain decomposition methods - algorithms and theory (R. Bank,R. L. Graham, J. Stoer, R. Varga, and H. Yserentant, eds.), Springer-Verlag, 2005

[31] L. N. Trefethen. Spectral Methods in MATLAB. Society for Industrial and Applied Mathematics, 2000

[32] Holger Wendland. Scattered data approximation. Cambridge University Press (Cambridge Monographs on Applied and Computational Mathematics, 17), 2005